\documentclass[12pt]{article}

\usepackage{array,multirow}
\usepackage{slashbox,pict2e}

\usepackage[text={6.25in,8.75in},centering]{geometry}

\usepackage{amsmath,amssymb}
\usepackage{graphicx}        
\usepackage{multicol}        
\usepackage{url}
 \newtheorem{theorem}{Theorem}[section]
 
 \newtheorem{proposition}[theorem]{Proposition}

 \newtheorem{algorithm}[theorem]{Algorithm}

\usepackage{bm}
\usepackage{enumerate}
\usepackage{dsfont}
\usepackage[algo2e,linesnumbered,ruled]{algorithm2e}
\usepackage{float}
\usepackage{perpage}
\MakeSorted{figure}
\MakeSorted{table}

\usepackage[pagewise]{lineno} 
\usepackage{verbatim}
\usepackage[config,labelfont={rm,rm},textfont=rm]{caption,subfig}
\usepackage{color}
\usepackage{wrapfig}

\renewcommand{\PackageWarningNoLine}[2]{}
\usepackage{amsmath}
\usepackage{amsfonts}
\usepackage{bm}
 
\newcommand{\av}[1]{\{#1\}} 

\newcommand{\jump}[1]{\lbrack\!\lbrack\,#1\,\rbrack\!\rbrack} 
\newcommand{\Eh}{{{\mathcal E}_h}} 
\newcommand{\Eho}{{{\mathcal E}^{o}_h}} 

\newcommand{\dyle}{\displaystyle}
\newcommand{\Ehb}{{{\mathcal E}^{\partial}_h}} 

\renewcommand{\lor }{\longrightarrow}

\newcommand{\jtau}{\jump{\taub}} 
\newcommand{\taub}{{\boldsymbol {\tau}}}

\numberwithin{equation}{section}
\numberwithin{figure}{section}
\numberwithin{table}{section}

\newcommand{\reals}[1]{\textrm{I\!{R}}^{#1}}
\newcommand{\K}{T} 
\renewcommand{\O}{\Omega} 
\renewcommand{\div}{{\rm div}} 

\newcommand{\n}{{\bf n}} 
\newcommand{\calB}{\mathcal{B}}

\newcommand{\calP}{\mathcal{P}}

\newcommand{\calA}{ \mathcal{A}}
\newcommand{\Th}{ \mathcal{T}_h}
\newcommand{\triplenorm}[1]{\left\vert\kern-0.9pt\left\vert\kern-0.9pt\left\vert #1 \right\vert\kern-0.9pt\right\vert\kern-0.9pt\right\vert}  
\newcommand{\calZ}{ \mathcal{Z}}

\graphicspath{{Picpdf/}{Piceps/}}
\makeatletter
\def\section{\@startsection {section}{1}{\z@}{-3.5ex plus -1ex minus -.2ex}{0.5em}{\large\bf}}
\def\subsection{\@startsection{subsection}{2}{\z@}{0.5ex}{0.5em}{\normalsize\bf}}
\def\subsubsection{\@startsection{subsubsection}{3}{\z@}{0.5ex}{-1em}{\normalsize\bf}}
\def\paragraph{\@startsection{paragraph}{4}{\z@}{0.5ex}{-1em}{\normalsize\bf}}
\def\subparagraph{\@startsection{subparagraph}{4}{\z@}{0.5ex}{-1em}{\normalsize\bf}}
\makeatother

\begin{document}

\title{A Block Solver for the Exponentially Fitted IIPG-0 method}
\author{Blanca Ayuso de Dios\footnotemark[1]\footnote{Centre de Recerca Matem\`atica,
   Barcelona, Spain. \texttt{bayuso@crm.cat}}
\and
Ariel Lombardi\footnotemark[2]\footnote{Departamento de Matem\'atica, Universidad de Buenos Aires \& CONICET,  Argentina. \texttt{aldoc7@dm.uba.ar}}
\and Paola Pietra
\footnotemark[3]\footnote{IMATI-CNR, Pavia, Italy,
  \texttt{pietra@imati.cnr.it}}
\and
Ludmil Zikatanov\footnotemark[4]\footnote{Department of Mathematics, Penn State
  University, USA \texttt{ltz@math.psu.edu}}
}
\date{}
%
%
\maketitle
\begin{abstract}
We consider an exponentially fitted discontinuous Galerkin method and
propose a robust  block solver for the resulting linear systems.
\end{abstract} 

\section{Introduction}
Let $\O\subset \reals{2}$ be a convex polygon, $f\in L^{2}(\O), g\in H^{1/2}(\partial\O)$ and let $\epsilon >0$ be constant. We consider the advection-diffusion problem
\begin{equation}\label{ayuso_mini_M8:1}
 -\div (\epsilon \nabla u - \beta u) = f
\quad \mbox{in } \, \Omega, \qquad
u=g \quad \mbox{on } \, \partial \O,
\end{equation}
where $\beta \in W^{1,\infty}(\Omega)$ derives from a potential
$\beta=\nabla \psi$. In applications to semiconductor devices, $u$
represents the concentration of positive charges, $\psi$ the
electrostatic potential and the electric field $|\nabla \psi|$ might
be fairly large in some parts of $\O$, so that \eqref{ayuso_mini_M8:1}
becomes advection dominated. Its robust numerical approximation and
the design of efficient solvers, are still nowdays a challenge.
Exponential fitting \cite{MR2143849} and discontinuous Galerkin (DG)
are two different approaches that have proved their usefulness for the
approximation of \eqref{ayuso_mini_M8:1}. Both methodologies have been
combined  in \cite{ariel-paola} to develop a new family of exponentially fitted DG methods (in primal and mixed formulation). In this note, we consider a variant of these schemes, based on the use of the Incomplete Interior Penalty IIPG-0 method and propose also an efficient block solver for the resulting linear systems.\\
\indent By introducing the change of variable 
\begin{equation}\label{ayuso_mini_M8:2}
\rho:= e^{-\frac{\psi}{\epsilon}} u 
\end{equation}
 problem \eqref{ayuso_mini_M8:1} can be rewritten as the following second order problem
\begin{equation}\label{ayuso_mini_M8:3}
-\nabla\cdot \left(\kappa \nabla \rho \right) =f  \mbox{   in  } \O, \quad \rho=\chi  \mbox{   on } \partial\O\;,
\end{equation}
where $\kappa:=\epsilon e^{\frac{\psi}{\epsilon}}$ and $\chi:=
e^{-\frac{\psi}{\epsilon}} g$. 
An IIPG-0 approximation to \eqref{ayuso_mini_M8:3} combined with a suitable
local approximation to  \eqref{ayuso_mini_M8:2}, 
gives rise to the EF-IIPG-0 scheme for  \eqref{ayuso_mini_M8:1}. 
We propose a block solver that uses ideas from \cite{MR2511726} to reduce the cost to that of a Crouziex-Raviart (CR) (exponentially fitted) discretization. By using Tarjan's algorithm, the associated matrix is further reduced to block lower triangular form, and a block Gauss-Siedel algorithm results in an exact solver. 

To give
a neat presentation, we focus on the case ${\bf \beta}=\nabla \psi$ piecewise
constant; $\psi$ piecewise linear continuous, although we
include some numerical results for a much more general case (cf. Test 2).
 Due to space restrictions, we describe the method and the solver and show some numerical results; further extensions of the method (allowing $\psi$ to be discontinuous) and the convergence analysis of the proposed solvers will be consider somewhere
else.
 \vskip -0.3cm
\section{The Exponentially Fitted IIPG-0 method}
\vskip -0.1cm
Let $\Th$ be a shape-regular
family of partitions of $\O$ into triangles $\K$ and let $h=\max_{\K \in \Th} h_{\K}$ with $h_{\K}$ denoting the diameter of $\K$ for each $\K \in \Th$. We assume $\Th$ does not contain hanging
nodes. 
We denote by $\Eho$ and $\Ehb$ the sets of all interior and boundary edges, respectively, and we set $\Eh=\Eho\cup \Ehb$.\\
\noindent {\bf  Average and jump trace operators:} Let $\K^{+}$ and $\K^{-}$ be two neighboring elements, and $\n^{+}$, $\n^{-}$ be their outward normal unit vectors, respectively ($\n^{\pm}=\n_{\K^{\pm}}$). Let
$\zeta^{\pm}$ and $\taub^{\pm}$ be the restriction of $\zeta$ and
$\taub$ to $\K^{\pm}$. We set:
\vskip -0.1cm 
\begin{equation*}
\begin{aligned}
2\av{\zeta}&=(\zeta^+ +\zeta^-),\quad
\jump{\zeta}=\zeta^+\n^++\zeta^-\n^- \quad&&\mbox{on } E\in
\Eho,\\
2\av{\taub}&=(\taub^+ +\taub^-), \quad
\jtau=\taub^+\cdot\n^+ + \taub^-\cdot\n^- &&\mbox{on } E\in
\Eho, \nonumber
\end{aligned}
\end{equation*}
and on  $e\in \Ehb$ we set $\jump{\zeta}=\zeta\n$ and $\av{\taub} = \taub$. We will also use the notation
\begin{equation*}
(u,w)_{\Th}=\dyle\sum_{\K \in \Th} \int_{\K} uw dx \qquad \langle u, w\rangle_{\Eh}=\dyle\sum_{e\in \Eh} \int_{e} u w ds \quad \forall\, u,w, \in V^{DG}\;,
\end{equation*}  
\noindent where  $V^{DG}$  is the discontinuous linear finite element space
defined by: 
\begin{equation*}
V^{DG}=\left\{  u\in L^{2}(\Omega) \,\, : \,\, u_{|_{\K}}\,\,  \in\,\,  \mathbb{P}^{1}(\K) \,\, \forall  \K \in \Th\,\right\},
\end{equation*}
$\mathbb{P}^{1}(\K)$ being the space of linear polynomials on $\K$. Similarly, $\mathbb{P}^{0}(\K)$ and $\mathbb{P}^{0}(e)$ are the spaces of constant polynomials on $\K$ and $e$, respectively.
For each $e\in \Eh$ (resp. for each $\K\in \Th$), let $\calP_{e}^{0}:L^{2}(e)\lor \mathbb{P}^{0}(e)$ (resp. $\calP_{\K}^{0}:L^{2}(\K)\lor \mathbb{P}^{0}(\K)$) be the $L^{2}$-orthogonal projection defined by  
\begin{equation*}
\calP_{e}^{0}(u):= \frac{1}{|e|} \int_{e} u, \quad \forall\, u\in L^{2}(e)\;, \quad\calP_{\K}^{0}(v):= \frac{1}{|\K|} \int_{\K} v, \quad \forall\, v\in L^{2}(\K)\;.
\end{equation*}
We denote by  $V^{CR}$ the classical Crouziex-Raviart (CR) space:
\begin{equation*}
V^{CR}\!=\!\left\{ v\in L^{2}(\Omega) \, : \, v_{|_{\K}}
\in\, \mathbb{P}^{1}(\K) \, \forall \K \in \Th\, \mbox{   and    }
\calP_{e}^{0} 
\jump{v}=0 \,\, \forall\, e\in \Eh\right\}.
\end{equation*}
Note that $v=0$ at the midpoint $m_{e}$ of each $e \in \Ehb$. To represent the functions in $V^{DG}$  we use the basis $\{\varphi_{e,\K}\}_{\K\in \Th, e\in \Eh}$, defined by
\begin{equation}\label{ayuso_mini_M8:4}
\forall \, \K\in \Th \quad \varphi_{e ,\K}(x) \in \mathbb{P}^{1}(\K)\quad e\subset \partial \K \quad \varphi_{e,\K}(m_{e'})=\delta_{e,e'} \quad \forall e' \in \Eh\;.
\end{equation}
In particular, any  $w\in  \mathbb{P}^{1}(\K)$ can be  written as $w=\sum_{e\subset \partial\K} w(m_e) \varphi_{e,\K}$.
\vskip -0.1cm
\subparagraph{The Exponentially fitted IIPG-0 method} 
We first consider the IIPG-0 approximation to the solution of \eqref{ayuso_mini_M8:3}:  Find $\rho \in V^{DG}$ such that $\calA(\rho,w)=(f,w)_{\Th}$  forall $w\in V^{DG}$  with 
\begin{equation}\label{ayuso_mini_M8:5}
\calA(\rho,w)= (\kappa^{\ast}_{\K} \nabla \rho,\nabla w)_{\Th} -\langle \av{ \kappa^{\ast}_{\K} \nabla \rho}, \jump{w}\rangle_{\Eh}+\langle S_{e} \av{ \jump{ \rho}}, \calP^{0}(\jump{w})\rangle_{\Eh}\;.
\end{equation}
Here, $S_e$ is the penalty parameter and $\kappa^{\ast}_{\K} \in \mathbb{P}^{0}(\K)$  the harmonic average approximation to $\kappa =\epsilon e^{\psi/\epsilon}$ both defined by \cite{ariel-paola}:
\begin{equation}\label{ayuso_mini_M8:6}
  \kappa^{\ast}_{\K}:=\frac{1}{ \mathcal{P}^{0}_{\K} (\kappa^{-1}) }= \frac{\epsilon }{  \mathcal{P}^{0}_{\K} (e^{-\frac{\psi}{\epsilon}} ) }\;, \qquad S_{e}:=  \alpha_e h_e^{-1}  \av{\kappa^{\ast}_{\K}}_{e} \;,
  \end{equation} 
Next, following \cite{ariel-paola} we introduce the local operator $\mathfrak{T} : V^{DG}\lor V^{DG}$ that approximates the change of variable \eqref{ayuso_mini_M8:2}:
\begin{equation}\label{ayuso_mini_M8:7}
\mathfrak{T}w :=\sum_{\K\in \Th} (\mathfrak{T} w)|_{\K}=\sum_{\K\in \Th} \sum_{e\subset \partial\K}
\calP^{0}_e(e^{-\frac{\psi}{\epsilon}})
w(m_e)\varphi_{e,\K} \quad \forall\, w\in V^{DG}\;.
\end{equation}
By setting  $\rho:= \mathfrak{T} u$ in \eqref{ayuso_mini_M8:5}, we finally get the EF-IIPG-0 approximation to~\eqref{ayuso_mini_M8:1}:\\   \indent    Find $u_h \in V^{DG}$ s.t. $\calB(u_h,w):=\calA(\mathfrak{T} u_h,w)=(f,w)_{\Th}$  $\forall\, w\in V^{DG}$ with 
\begin{equation}\label{ayuso_mini_M8:8}
\calB(u,w)\!=\!\!(\kappa^{\ast}_{\K} \nabla \mathfrak{T} u,\nabla w)_{\Th}\! -\langle\av{ \kappa^{\ast}_{\K} \nabla \mathfrak{T} u}, \jump{w}\rangle_{\Eh}\!+ \langle S_{e} \av{ \jump{ \mathfrak{T} u}}, \calP^{0}\jump{w}\rangle_{\Eh}\;.
\end{equation}
It is important to emphasize that the use of harmonic average to approximate $\kappa=\epsilon e^{\psi/\epsilon}$ as defined in \eqref{ayuso_mini_M8:6}  together with the definition of the local approximation of the change of variables prevents possible overflows in the computations when $\psi$ is large and $\epsilon$ is small. (See \cite{ariel-paola} for further discussion). 
Also, these two ingredients are essential to ensure that the resulting method  has an automatic upwind mechanism built-in that allows for an accurate approximation of the solution of \eqref{ayuso_mini_M8:1} in the advection dominated regime. We will discuss this in more detail in Section \ref{ayuso_mini_M8:00}.

Prior to close this section, we define for each $e\in \Eh$ and $\K\in \Th$:
\begin{equation*}
 \psi_{m,e}:=\min_{x\in e} \psi(x) \quad  \psi_{m,\K}:= \min_{x\in \K} \psi(x); \quad\, \psi_{m,\K}\leq \psi_{m,e} \mbox{   for   }  e\subset \partial\K\;.
 \end{equation*}
In the advection dominated regime $ \epsilon \ll |\beta| h= |\nabla \psi| h$  
 \begin{equation}\label{ayuso_mini_M8:9}
 \calP_{\K}^{0} (e^{-(\psi/\epsilon)}) \simeq \epsilon^{2} e^{-\frac{\psi_{m,\K}}{\epsilon}} \qquad\qquad  \calP^{0}_{e_i} (e^{-\psi/\epsilon}) \simeq \epsilon \, e^{-\frac{\psi_{m,e}}{\epsilon}}\;.
  \end{equation}
 The first of the above scalings together with the definitions in \eqref{ayuso_mini_M8:6} implies
 \begin{equation}\label{ayuso_mini_M8:01}
 \kappa_{\K}^{\ast} \simeq \frac{1}{\epsilon}e^{\frac{\psi_{m,\K}}{\epsilon}}\;,  \qquad S_e \simeq \frac{\alpha}{2\epsilon}|e|^{-1} e^{\frac{(\psi_{m,\K_1}+\psi_{m,\K_2})}{\epsilon}}\; \quad e=\partial\K_1\cap \partial\K_2\;.
 \end{equation}
 \vglue -0.2cm
\vskip -0.3cm
\section{Algebraic System \& Properties }
\label{ayuso_mini_M8:00}
\vskip -0.1cm
Let $A$ and $B$ be the operators associated to the bilinear forms $\calA(\cdot,\cdot)$ \eqref{ayuso_mini_M8:5} and $\calB(\cdot,\cdot)$  \eqref{ayuso_mini_M8:8}, respectively. We denote by $\mathbb{A}$ and $\mathbb{B}$ their matrix representation in the basis  $\{\varphi_{e,\K}\}_{\K\in \Th, e\in \Eh}$  \eqref{ayuso_mini_M8:4}.  In this basis, the operator $\mathfrak{T}$ defined in \eqref{ayuso_mini_M8:7} is represented as a diagonal matrix, $\mathbb{D}$, and  $\mathbb{B}=\mathbb{A}\mathbb{D}$. Thus, the approximation to \eqref{ayuso_mini_M8:3} and \eqref{ayuso_mini_M8:1} amounts to solve the linear systems (of dimension $2n_e-n_b$; with $n_e$ and $n_b$ being the cardinality of $\Eh$ and $\Ehb$, respectively):
\begin{equation}\label{ayuso_mini_M8:02}
\mathbb{A} \bm{\rho}=\bm{F}\;, \quad\mbox{  and   } \quad \mathbb{D}\bm{u}=\bm{\rho}  \qquad\mbox{ or  }\qquad \mathbb{B} \bm{u}= \widetilde{\bm{F}}\;, 
\end{equation}
where $\bm{\rho}, \bm{u}, \bm{F}$ and $\widetilde{\bm{F}} $ are the vector representations of  $\rho, u$ and the rhs of the approximate problems.
From the definition \eqref{ayuso_mini_M8:7} of $\mathfrak{T}$ it is easy to deduce the scaling of the entries of the diagonal matrix $\mathbb{D}=(d_{i,i})_{i=1}^{2n_e-n_b}$.
 \begin{equation*}
 \mathbb{D}=(d_{i,j})_{i,j=1}^{2n_e-n_b}  \quad d_{i,i}= \calP^{0}_{e_i} (e^{-\psi/\epsilon}) \simeq \epsilon \, e^{-\frac{\psi_{m,e}}{\epsilon}}\;, \quad d_{i,j}\equiv 0 \quad i\ne j\;.
 \end{equation*}
 We now revise a result from \cite{MR2511726}:
\begin{proposition}\label{ayuso_mini_M8:03}
Let $\calZ\subset V^{DG}$ be the space defined by 
\begin{equation*}
\calZ=\left\{ z\in L^{2}(\Omega) \, : \, z_{|_{\K}}\,
\in\, \mathbb{P}^{1}(\K) \,\, \forall \K \in \Th\, \mbox{  and  }
\calP_{e}^{0} 
\av{v}=0 \,\, \forall\, e\in \Eho \right\}.
\end{equation*}
Then, for any $w\in V^{DG}$ there exists a unique $w^{cr}\in V^{CR}$ and a unique $w^{z}\in \calZ$ such that $w=w^{cr}+w^{z}$ , that is:  $V^{DG}=V^{CR}\oplus \calZ$. Moreover,
$ \calA(w^{cr},w^{z})=0$  $\forall\, w^{cr} \in V^{CR},$ and  $\forall\, w^{z} \in \calZ$.
\end{proposition}
Proposition \ref{ayuso_mini_M8:03} provides a simple {\it change of  basis}  from  $\{\varphi_{e ,\K}\}$ to canonical basis in $V^{CR}$ and $\calZ$ that results in the following algebraic structure for \eqref{ayuso_mini_M8:02}:
\begin{equation}\label{ayuso_mini_M8:04}
\bm{\rho}=\left[\begin{array}{cccc}
&\bm{\rho^{z}}\\
&\bm{\rho^{cr}}\end{array}\right], \qquad
\mathbb{A}= \left[ \begin{array}{cccc}
&\mathbb{A}^{zz} & \bm{0} &\\
& \mathbb{A}^{vz} & \mathbb{A}^{vv} &
\end{array}\right], \qquad \mathbb{B}= \left[ \begin{array}{cccc}
&\mathbb{B}^{zz} &  0 &\\
& \mathbb{B}^{vz} & \mathbb{B}^{vv} &
\end{array}\right]. 
\end{equation}
Due to the assumed continuity of $\psi$,  $\mathbb{D}$ is still diagonal in this basis. The algebraic structure \eqref{ayuso_mini_M8:04} suggests the following exact solver:  
\begin{algorithm}
Let $u_0$ be a given initial guess. For $k\ge 0$, and given
$u_k=z_k+v_k$, the next iterate
$u_{k+1}=z_{k+1}+v_{k+1}$ is defined via the two steps:
\begin{enumerate}
\item[1.] Solve $\calB(u^{z}_{k+1},w^{z}) = (f,w^{z})_{\Th}
 \quad \forall\, w^{z}\in \calZ$.
\item[2.] Solve $\calB(u^{cr}_{k+1},w^{cr}) = (f,w^{cr})_{\Th}-\calB(u^{z}_{k+1},w^{cr})
\quad \forall\, w^{cr}\in V^{CR}$.
\end{enumerate}
\end{algorithm}
Next, wet discuss how to solve efficiently each of the above steps:\\

\noindent {\bf Step 2: Solution in $V^{CR}$.}  In \cite{MR2511726} it was shown that the block  $\mathbb{A}^{vv}$ coincides with the stiffness matrix of a CR discretization of \eqref{ayuso_mini_M8:3}, and so it is an s.p.d. matrix. However, this is no longer true for $\mathbb{B}^{vv}$ which is positive definite but non-symmetric.
\begin{equation*}
\calB(u^{cr},w^{cr})=(\kappa^{\ast}_{\K}\nabla \mathfrak{T} u^{cr},\nabla w^{cr})_{\Th} \quad \forall\, \,u^{cr}\, ,w^{cr} \in V^{CR}\;.
\end{equation*}
In principle, the sparsity pattern of $\mathbb{B}^{vv}$ is that of a symmetric matrix. Using \eqref{ayuso_mini_M8:9} and \eqref{ayuso_mini_M8:6}, we find that the entries of  the matrix scale as:
\begin{equation}\label{ayuso_mini_M8:05}
  \mathbb{B}^{vv}=\left(b^{cr}_{i,j}\right)_{i,j}^{n_{cr}:=n_e-n_b} \quad b^{cr}_{i,j}:= \kappa_{\K}^{\ast}\frac{|e_i| |e_j|}{|\K|} \,{\bf n}_{e_i} \cdot {\bf n}_{e_j} d_{j} \simeq e^{-\frac{(\psi_{m,e}-\psi_{m,\K})}{\epsilon}} 
 \end{equation}
Since $\psi$ is assumed to be piecewise linear, for each $\K$, it attains its minimum (and also its maximum) at a vertex of $\K$, say $\bm{x_0}$ and $\psi_{m,e}$ is attained at one of the vertex of the edge $e$, say $\bm{x_e}$. In particular, this implies that
\begin{equation*}
 \psi_{m,e}-\psi_{m,\K} \approx   \nabla \psi \cdot (\bm{x_e}-\bm{x_0}) =\beta \cdot (\bm{x_e}-\bm{x_0})=\left\{ \begin{array}{ccc}
 &0 &\quad  \bm{x_e}=\bm{x_0}\\
  &|\beta| h  &\quad  \bm{x_e}\ne \bm{x_0}
  \end{array}\right.
 \end{equation*}
Hence, in the advection dominated case $\epsilon\ll |\beta| h$ some of the entries in \eqref{ayuso_mini_M8:05}  vanish (up to machine precision) for $\epsilon$ small; this is  the automatic upwind mechanism intrinsic of the method. As a consequence, the sparsity pattern of $\mathbb{B}^{vv}$ is no longer symmetric and this can be exploited to re-order the unknowns so that $\mathbb{B}^{vv}$ can be reduced to block lower triangular form.\\
 Notice also that for $\Th$ acute, the block $\mathbb{A}^{vv}$  being the stiffness matrix
 of the Crouziex-Raviart approximation to \eqref{ayuso_mini_M8:3}, is an M-matrix. Hence, since the
 block $\mathbb{B}^{vv}$ is the product of a positive diagonal matrix
 and $\mathbb{A}^{vv}$, it will also be an $M$-matrix if
 the triangulation is acute (see  \cite{MR2143849}).\\

\noindent {\bf Step 1: Solution in the $\calZ$-space.} In \cite{MR2511726} it was
 shown that $A^{zz}$ is a diagonal p.d. 
 matrix. This is also true for  $\mathbb{B}^{zz}$ since it is the product of two diagonal matrices. The continuity of $\psi$ implies
\begin{equation}
\calB(u^{z},w^{z})=\langle S_{e}\mathfrak{T}\jump{u^{z}}, \calP^{0}_{e} (\jump{w^{z}}) \rangle_{\Eh}  \quad \forall\,\, u^{z}, \, w^{z} \in \calZ\;.
\end{equation}
Using \eqref{ayuso_mini_M8:9} and \eqref{ayuso_mini_M8:6} we observe that the entries of $\mathbb{B}^{zz}$ scale as:
\begin{equation*}
\mathbb{B}^{zz}= \left( b_{i,j}\right)_{i=1}^{n_e} \quad b_{i,j}= S_{e_i} |e_{i}| d_{j}\delta_{i,j}   \simeq \delta_{i,j}\frac{\alpha}{2}\, e^{-(\psi_{m,e}- \psi_{m,\K_1}-\psi_{m,\K_2})/\epsilon}
\end{equation*}
which are always positive, so in particular $\mathbb{B}^{zz}$ it is also an $M$-matrix.

\section{Block Gauss-Siedel solver for $V^{CR}$-block}

We now consider re-orderings of the unknowns (dofs),
which reduce $\mathbb{B}^{vv}$ to block lower triangular form. For
such reduction, we use the algorithm from~\cite{MR0304178} which
roughly amounts to \emph{partitioning} the set of dofs into
non-overlapping blocks.  In the strongly advection dominated case the
size of the resulting blocks is small and a block Gauss-Seidel method
is an efficient solver.  Such techniques have been studied 
in~\cite{MR1718671} for conforming methods.  The idea is to consider
the \emph{directed} graph $\bm{G}=(\bm{V},\bm{E})$ associated with
$\mathbb{B}^{vv} \in \reals{n_{cr}\times n_{cr}}$; $\bm{G}$ has
$n_{cr}$ vertices labeled $\bm{V}=\{1,\ldots,n_{cr}\}$ and its set of
edges \emph{edges} $\bm{E}$ has cardinality equal to the number of
nonzero entries\footnote{Each dof corresponds to a vertex in the
  graph; each nonzero entry to an edge.} of $\mathbb{B}^{vv}$.  By
definition, $(i,j)\in \bm{E}$ \emph{iff} $b^{cr}_{ij}\neq 0$. Note
that in the advection dominated case, due to the nonsymmetric pattern
of $\mathbb{B}^{vv}$ (caused by the built-in upwind mechanism), we may
have $(i,j)\in \bm{E}$, while $(j,i)\notin \bm{E}$. Then, the problem
of reducing $\mathbb{B}^{vv}$ to block lower triangular form of
$\mathbb{B}^{vv}$ is equivalent to partitioning $\bm{G}$ as a union of
strongly connected components. Such partitioning induces
non-overlapping partitioning of the set of dofs,
$\bm{V}=\cup_{i=1}^{N_b}\omega_i$. For $i=1,\ldots, N_b$, let $m_i$
denote the cardinality of $\omega_i$; let $\mathbb{I}_i\in
\reals{n_{cr}\times m_i}$ be the matrix that is identity on dofs
in $\omega_i$ and zero otherwise; and
$\mathbb{B}^{vv}_i=\mathbb{I}_i^T\mathbb{B}^{vv}\mathbb{I}_i$ is the
 block corresponding to the dofs in $\omega_i$.  The
block Gauss-Seidel algorithm reads: \emph{Let $\bm{u}^{cr}_0$ be
  given, and assume $\bm{u}^{cr}_{k}$ has been obtained. Then
  $\bm{u}^{cr}_{k+1}$ is computed via:} For $i=1,\ldots N_b$
 \begin{equation}\label{ayuso_mini_M8:06}
 \bm{u}^{cr}_{k+i/N_b}=\bm{u}^{cr}_{k+(i-1)/N_b} +\mathbb{I}_i(\mathbb{B}^{vv}_i)^{-1}\mathbb{I}_i^T\left(\bm{F}-\mathbb{B}^{vv}\bm{u}^{cr}_{k+(i-1)/N_b}\right)\;.
\end{equation}
As we report in Section \ref{ayuso_mini_M8:07}, in the advection
dominated regime the action of $(\mathbb{B}^{vv}_i)^{-1}$ can be
computed exactly since the size of the blocks $\mathbb{B}^{vv}_i$ is
small.
\section{Numerical Results}\label{ayuso_mini_M8:07}

\begin{figure}[!htb]
\centering
  \subfloat[Test~1 with $\epsilon=10^{-5}$]
  {\includegraphics*[width=0.4\textwidth]{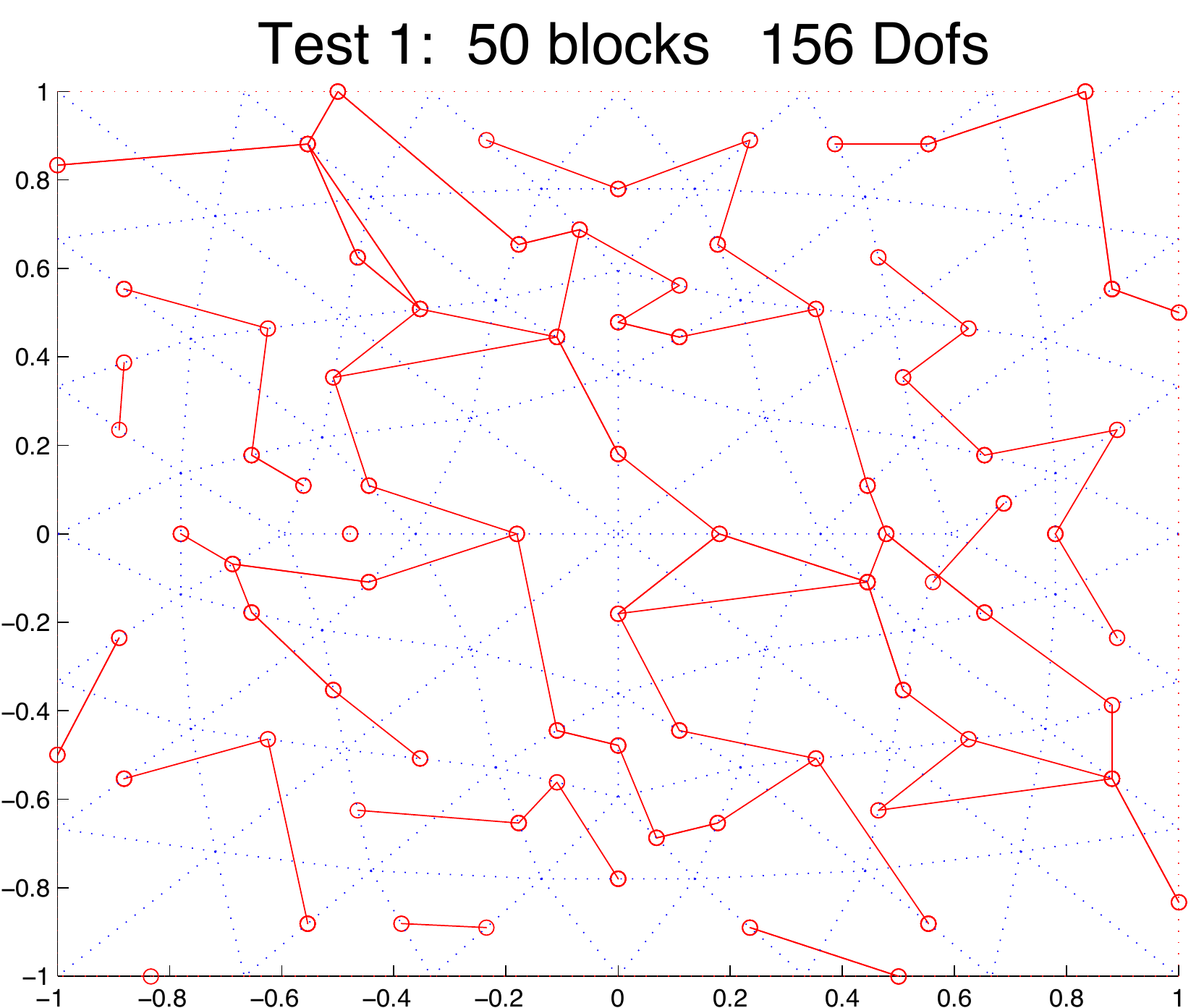}}  
\hspace*{0.1\textwidth}
  \subfloat[Test~2 with  $\epsilon=10^{-7}$]
  {\includegraphics*[width=0.4\textwidth]{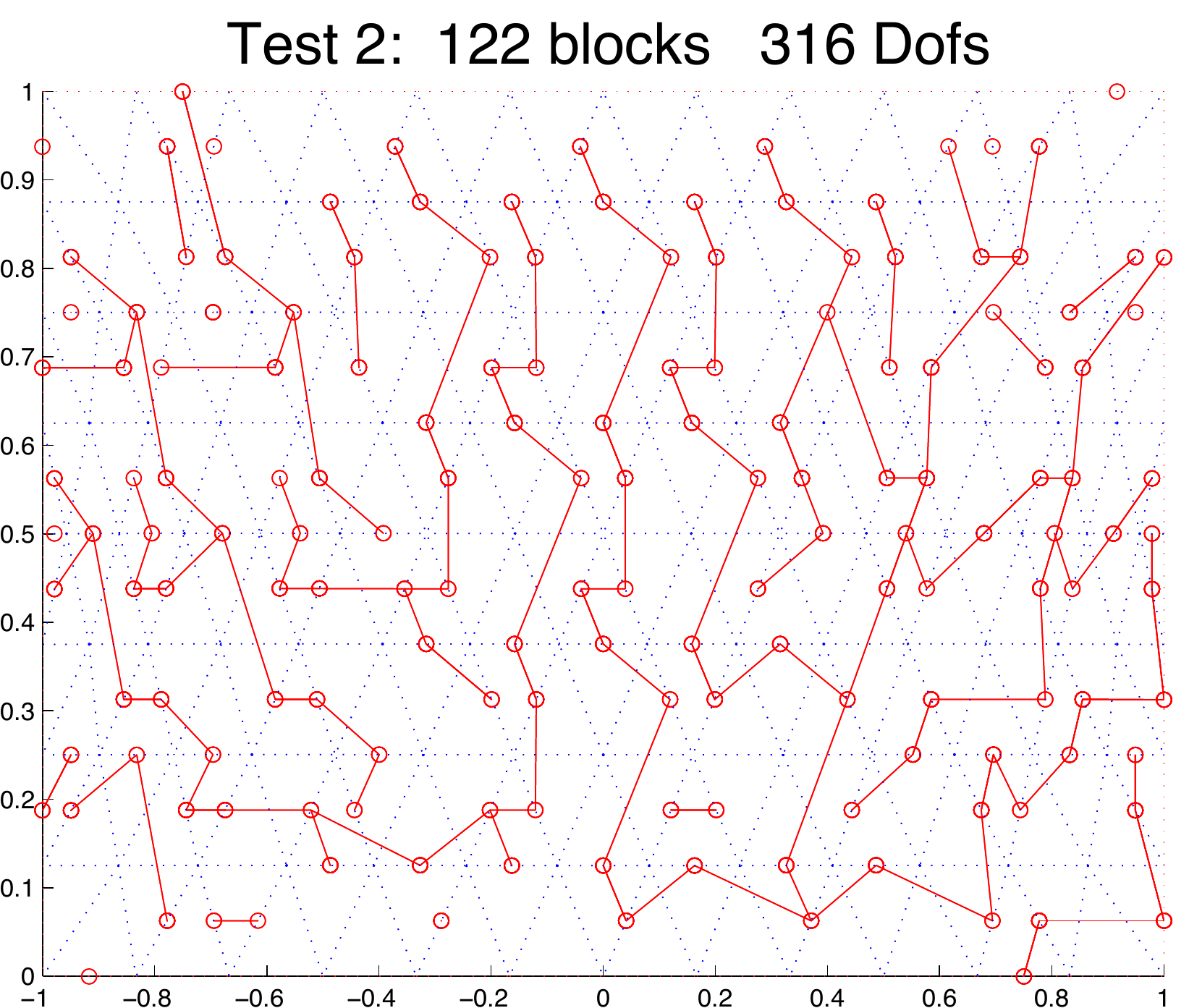}}
  \caption{Plot of the connected components
    (blocks) of $\mathbb{B}^{vv}$ created during
    Tarjan's algorithm. \label{ayuso_mini_M8:f0}}
\end{figure}
 We present a set of numerical experiments to assess the
performance of the proposed block solver. The tests refer to problem \eqref{ayuso_mini_M8:3}  with $\epsilon=10^{-3},10^{-5},10^{-7}$,  and $\O$ is  triangulated
with a family of unstructured triangulations $\Th$. In the tables
given below $J=1$ corresponds to the coarsest grid and each refined
triangulation on level $J$, $J=2,3,4$ is obtained by subdividing
each of the $\K\in \Th$ on level $(J-1)$ into four
congruent triangles. From the number of triangles $n_{\K}$ the total number of dofs for the DG approximation is $3n_{\K}$ and $n_{e}-n_{b}$ for the CR part of the solution.\\


\noindent {\bf Test 1. Boundary Layer:}  $\Omega=(-1,1)^2$, $\beta=[1,1]^{t}$, $n_{\K}=112$ for the coarsest mesh and $f$ is such that the exact solution is given by
\begin{equation*}
u(x,y)=\left( x+\frac{1+e^{-2/\epsilon} -2e^{(x-1)/\epsilon}}{1-e^{-2/\epsilon}}\right)\left( y+\frac{1+e^{-2/\epsilon} -2e^{(y-1)/\epsilon}}{1-e^{-2/\epsilon}}\right)\;.
\end{equation*} 
\noindent {\bf Test 2. Rotating Flow:}  $\Omega=(-1,1)^\times (0,1)$, $f=0$ and curl$\beta\ne 0$,
\begin{equation*}
\beta=\left[
\begin{array}{r}
2y(1-x^{2})\\
-2x(1-y^{2})\end{array}\right]^{t} \,\,\, g(x,y)=\left\{\begin{array}{llll}
1+\tanh{(10(2x+1))} & \,\, x\leq 0,\,\, y=0, \\
0 &   \mbox{  elsewhere}\;.
\end{array}\right.
\end{equation*}
We stress that this test does not fit in the simple description given here, and special care is required (see \cite{ariel-paola}). For the approximation, for each $\K\in \Th$, with barycenter $(x_{\K},y_{\K})$,  we use the approximation $\beta|_{\K} \approx \nabla \psi|_{\K}$ with $\psi|_{\K}=2y_{\K}(1-x^{2}_{\K})x-2x_{\K}(1-2y_{\K}^{2})y$ (and so $\psi$ discontinuous). The coarsest grid has $n_{\K}=224$ triangles.

In Figure~\ref{ayuso_mini_M8:f0} the plot of the
connected components of the graph depicting the blocks for
$\mathbb{B}^{vv}$ created during Tarjan's algorithm, on the coarsest
meshes is shown; for Test 1 with $\epsilon=10^{-5}$
and for Test 2 with $\epsilon=10^{-7}$.  In Tables
\ref{ayuso_mini_M8:t0} 
are given, the number of blocks $N_{b}$ created
during Tarjan's algorithm. 
We also report in this table the size of
the largest block created ($M_{b}$ maximum size) and the average size
of the blocks $n_{av}$. Observe that in the advection dominated regime
the largest block has a very small size compared to the total size of
the system.  After Tarjan's algorithm is used to re-order the matrix
$\mathbb{B}^{vv}$, we use the block Gauss-Seidel algorithm
\eqref{ayuso_mini_M8:06} where each small block is solved exactly.

\begin{table}
\centering
\subfloat[Test 1]{
\scriptsize
\begin{tabular}{|c|c||c|c|c|c|}
\hline\hline
\backslashbox{$\epsilon$}{$J$} & 
&\phantom{xx} $1$ \phantom{xx}  &\phantom{xx} $2$\phantom{xx}    & \phantom{xx}$3$\phantom{xx}   &\phantom{xx} $4$\phantom{xxx}  \\
\hline\hline
\multirow{3}{*}{$10^{-3}$}
 & $N_b$ &  44 & 150 &  484 & 1182\\ 
 & $M_b$ & 23   &  47    &     95  &      191  \\
 & $n_{av}$ & 3.55 &  4.32 &    5.45 &   9.02 \\\hline
\multirow{3}{*}{$10^{-5}$}
 & $N_b$  &  50 &    210  &  866 &  3474 \\  
 & $M_b$ & 23   &  47    &     95  &      191  \\
 & $n_{av}$ & 3.12 &   3.08 &   3.05 &    3.07  \\\hline
\multirow{3}{*}{$\!\!\!10^{-7}\!\!\!$}
  & $N_b$ & 50   &   210  &866 & 3522   \\  
  & $M_b$ & 23   &  47    &     95  &      191  \\
 & $n_{av}$ & 3.12 &  3.08 &   3.05 &    3.03   \\\hline
\hline
\end{tabular}
}
\subfloat[Test 2]{
\scriptsize
\begin{tabular}{|c|c||c|c|c|c|}
\hline\hline
\backslashbox{$\epsilon$}{$J$} & 
 &\phantom{xx} $1$ \phantom{xx}  & \phantom{xx}$2$\phantom{xx}    & \phantom{xx}$3$\phantom{xx}   &\phantom{xx} $4$\phantom{xxx} \\
\hline\hline
\multirow{3}{*}{$10^{-3}$}
 & $N_b$ &  31  &  1   &   1  &  1\\ 
  & $M_b$ & 211 & 1304 & 5296  & 21344 \\
 & $n_{av}$ & 10.19 &  1304 &  5296 &  21344 \\\hline
\multirow{3}{*}{$10^{-5}$}
 & $N_b$  &  122  &  468 &1822 & 7106 \\
 & $M_b$ & 4    &4    & 7  &37 \\
 & $n_{av}$ & 2.59 & 2.78 &    2.91 &    3.00  \\\hline
\multirow{3}{*}{$10^{-7}$}
  & $N_b$ & 122  &  468   & 1832 &   7247   \\
   & $M_b$ &4  &  4 &       4 &          6 \\
 & $n_{av}$ & 2.59 &    2.78 &    2.89 &    2.95    \\\hline
\hline
\end{tabular}
}
\caption{Number of blocks ($N_b)$ created during the Tarjan's ordering
  algorithm, size of largest block ($M_b$) and  average size of blocks $(n_{av})$.\label{ayuso_mini_M8:t0}}
\end{table} 
\section*{Acknowledgments}
This work started while the first two authors were visiting the
IMATI-CNR, Pavia in October 2010. Thanks go to the IMATI for the
hospitality and support. The first author was partially supported by
MEC grant MTM2008-03541, the second author was supported by CONICET
and the fourth author is supported in part by the National Science
Foundation NSF-DMS 0810982.




\end{document}